\documentclass{gtart_h}

\def\ifplaintex{\expandafter\ifx\csname documentclass\endcsname\relax}

\def\gtp{{\mathsurround=0pt\it $\cal G\mskip-2mu$eometry \&\ 
$\cal T\!\!$opology $\cal P\!$ublications}}  

\def\recd{{\small Received:\qua\receiveddate\ifx\reviseddate\relax
\else\qquad Revised:\qua\reviseddate\fi\par}} 

\def\lognumber#1{\def\thelognumber{#1}}
\def\volumenumber#1{\def\thevolumenumber{#1}}
\def\volumeyear#1{\def\thevolumeyear{#1}}
\def\papernumber#1{\def\thepapernumber{#1}}
\def\pagenumbers#1#2{\def\startpage{#1}\def\finishpage{#2}}
\def\published#1{\def\publishdate{#1}}

\def\received#1{\def\receiveddate{#1}}
\def\revised#1{\def\reviseddate{#1}}
\def\accepted#1{\def\accepteddate{#1}}
\def\asciititle#1{\def\theasciititle{#1}}
\def\covertitle#1{\def\thecovertitle{#1}}

\long\def\asciiabstract#1{\long\def\theasciiabstract{#1}}

\let\\\par\let\thelognumber\relax\let\thevolumenumber\relax
\let\thepapernumber\relax\let\thevolumeyear\relax\let\startpage\relax
\let\finishpage\relax\let\publishdate\relax\let\receiveddate\relax
\let\reviseddate\relax\let\accepteddate\relax\let\theasciititle\relax
\let\thecovertitle\relax\let\theasciiauthors\relax
\let\theasciiabstract\relax

\let\theasciiemail\relax

\ifplaintex
\font\logobig=cmssbx10 scaled 3836
\font\logomed=cmssbx10 scaled 2557
\else
\font\logobig=cmssbx10 scaled 4200
\font\logomed=cmssbx10 scaled 2800
\fi

\long\def\makeagttitle{   
\count0=\startpage
\agt\hfill      
\hbox to 45truept{\vbox to 0pt{\vglue -13truept{\logomed A\kern -.37em{\logobig 
T}\kern -.38em G}\vss}\hss}
\break
{\small Volume \thevolumenumber\ (\thevolumeyear)
\startpage--\finishpage\nl
Published: \publishdate}

\vglue .25truein

{\parskip=0pt\leftskip 0pt plus
1fil\def\\{\par\smallskip}{\Large\bf\thetitle}\par\medskip} \vglue
0.05truein

{\parskip=0pt\leftskip 0pt plus 1fil\def\\{\par}{\sc\theauthors}
\par\medskip}%
 
\vglue 0.03truein 

{\small\leftskip 25truept\rightskip 25truept{\bf Abstract}\stdspace\theabstract

{\bf AMS Classification}\stdspace\theprimaryclass
\ifx\thesecondaryclass\relax\else; \thesecondaryclass\fi\par
{\bf Keywords}\stdspace \thekeywords\par}\vglue 7truept

}   

\ifplaintex
\font\phead=cmsl9 scaled 950
\font\pnum=cmbx10 scaled 913
\font\pfoot=cmsl9 scaled 950
\headline{\vbox to 0pt{\vskip -4.5mm\line{\small\phead\ifnum
\count0=\startpage ISSN 1472-2739 (on-line) 1472-2747 (printed)
\hfill {\pnum\folio}\else\ifodd\count0\def\\{ }%
\ifx\theshorttitle\relax\thetitle\else\theshorttitle\fi\hfill{\pnum\folio}
\else\def\\{ and }{\pnum\folio}\hfill\ifx\theshortauthors\relax\theauthors
\else\theshortauthors\fi\fi\fi}\vss}}
\footline{\vbox to 0pt{\vglue 0mm\line{\small\pfoot\ifnum\count0=\startpage
\copyright\ \gtp\hfill\else
\agt, Volume \thevolumenumber\ (\thevolumeyear)\hfill\fi}\vss}}
\else
\font=cmsl9 scaled 1050
\font\lnum=cmbx10 
\font=cmsl9 scaled 1050
\makeatletter
\def\@oddhead{{\small\lhead\ifnum\count0=\startpage ISSN 1472-2739 
(on-line) 1472-2747 (printed)\hfill {\lnum\number\count0}\else\ifodd\count0
\def\\{ }\ifx\theshorttitle\relax \thetitle \else\theshorttitle\fi\hfill
{\lnum\number\count0}\else\def\\{ and }{\lnum\number\count0}
\hfill\ifx\theshortauthors\relax 
\theauthors\else\theshortauthors\fi\fi\fi}}\def\@evenhead{\@oddhead}
\def\@oddfoot{\small\lfoot\ifnum\count0=\startpage\copyright\ \gtp\hfill\else
\agt, Volume \thevolumenumber\ (\thevolumeyear)\hfill\fi}
\def\@evenfoot{\@oddfoot}
\makeatother
\fi
\let\maketitlepage\makeagttitle
\let\makeshorttitle\maketitlepage
\let\maketitle\maketitlepage

\newwrite\gtoutfile
\long\gdef\makeheadfile{  
{\def\\{, }\def\s{ }
\immediate\openout\gtoutfile head.xxx
\immediate\write\gtoutfile{Proxy-for: \ifx\theasciiauthors\relax
\theauthors\else\theasciiauthors\fi\s<\ifx\theasciiemail\relax\theemail\else\theasciiemail\fi>}
\immediate\write\gtoutfile{\noexpand\\}
\immediate\write\gtoutfile{Authors: \ifx\theasciiauthors\relax
\theauthors\else\theasciiauthors\fi}
{\def\\{ }\immediate\write\gtoutfile{Title: \ifx\theasciititle\relax
\thetitle\else\theasciititle\fi}}
\immediate\write\gtoutfile{Subj-class: GT or SG, GR etc}
\immediate\write\gtoutfile{MSC-class: \theprimaryclass\ifx\thesecondaryclass\relax\else, \thesecondaryclass\fi}
\immediate\write\gtoutfile{Journal-ref: Algebr. Geom. Topol. \thevolumenumber\s
(\thevolumeyear) \startpage-\finishpage}
\immediate\write\gtoutfile{Comments: Published by Algebraic and
Geometric Topology at}
\immediate\write\gtoutfile{\s\s\s  http://www.maths.warwick.ac.uk/agt/AGTVol\thevolumenumber/agt-\thevolumenumber-\thepapernumber.abs.html}
\immediate\write\gtoutfile{\noexpand\\}
\immediate\write\gtoutfile{}
\ifx\theasciiabstract\relax
\immediate\write\gtoutfile{\theabstract}\else
\immediate\write\gtoutfile{\theasciiabstract}\fi
\immediate\write\gtoutfile{}
\immediate\write\gtoutfile{\noexpand\\}
\immediate\write\gtoutfile{}
\immediate\closeout\gtoutfile}}  

\def\maketitlepage{\makeagttitle\makeheadfile}
\let\makeshorttitle\maketitlepage
\let\maketitle\maketitlepage

\lognumber{30}
\volumenumber{5}
\volumeyear{2005}
\papernumber{30}
\pagenumbers{725}{740}
\received{13 January 2005} 
\revised{5 May 2005}
\accepted{21 June 2005}
\published{13 July 2005}

\usepackage{amsmath}
\usepackage{amssymb}

\usepackage[all]{xy}

\theoremstyle{plain} 

\newtheorem{thm}{Theorem}[section]
\newtheorem{cor}[thm]{Corollary}

\newtheorem{lem}{Lemma}[section]
\newtheorem*{thm*}{Theorem}
\theoremstyle{definition}

\newcommand{\ra}{\ensuremath{\longrightarrow}}
\newcommand{\Q}{\ensuremath{\mathbb{Q}}}
\newcommand{\Z}{\ensuremath{\mathbb{Z}}}
\newcommand{\IA}{{\rm IA}}
\newcommand{\OA}{{\rm OA}}

\def\Wedge{\text{\Large\mathsurround0pt$\wedge$}}

\begin{document}

\title{The Johnson homomorphism and\\the second cohomology of $\IA_n$}
\covertitle{The Johnson homomorphism and\\the second cohomology of IA$_n$}
\asciititle{The Johnson homomorphism and\\the second cohomology of IA_n}
\author{Alexandra Pettet}
\address{Department of Mathematics, University of Chicago\\Chicago, IL 60637, USA}
\email{alexandra@math.uchicago.edu}

\begin{abstract} 
Let $F_n$ be the free group on $n$ generators. Define $\IA_n$ to be 
group of automorphisms of $F_n$ that act trivially on first homology. 
The Johnson homomorphism in this setting is a map from $\IA_n$ to 
its abelianization. The first goal of this paper is 
to determine how much this map contributes to 
the second rational cohomology of $\IA_n$. 

A descending central series of $\IA_n$ is given by the subgroups
$K_n^{(i)}$ which act trivially on $F_n/F_n^{(i+1)}$, the free rank
$n$, degree $i$ nilpotent group. It is a conjecture of Andreadakis
that $K_n^{(i)}$ is equal to the lower central series of $\IA_n$;
indeed $K_n^{(2)}$ is known to be the commutator subgroup of
$\IA_n$. We prove that the quotient group $K_n^{(3)}/\IA_n^{(3)}$ is
finite for all $n$ and trivial for $n=3$.  We also compute the rank of
$K_n^{(2)}/K_n^{(3)}$.
\end{abstract}
\asciiabstract{%
Let F_n be the free group on n generators. Define IA_n to be group
of automorphisms of F_n that act trivially on first homology.  The
Johnson homomorphism in this setting is a map from IA_n to its
abelianization. The first goal of this paper is to determine how much
this map contributes to the second rational cohomology of IA_n.  A
descending central series of IA_n is given by the subgroups K_n^(i)
which act trivially on F_n/F_n^(i+1), the free rank n, degree i
nilpotent group.  It is a conjecture of Andreadakis that K_n^(i) is
equal to the lower central series of IA_n; indeed K_n^(2) is known to
be the commutator subgroup of IA_n. We prove that the quotient group
K_n^(3)/IA_n^(3) is finite for all n and trivial for n=3.  We also
compute the rank of K_n^(2)/K_n^(3).}


\maketitle


\section{Introduction}
Let ${\rm Aut}(F_n)$ denote the automorphism group of the free group
$F_n$ on generators $\{x_1,\ldots,x_n\}$. 
The IA-automorphism group of $F_n$, denoted by $\IA_n$, 
is defined to be the kernel of the exact sequence 
\[
1 \ra {\rm IA}_n \ra {\rm Aut}(F_n) \ra 
{\rm GL}(n,\mathbb{Z}) \ra 1 
\]
induced by the action of ${\rm Aut}(F_n)$ on the abelianization of $F_n$. 

Although Nielsen studied $\IA_n$ as early as 1918, it is still the case that
little is understood about this group.  For 
instance, while there is a well-known finite set of generators, discovered by
Magnus in 1934 \cite{magnus}, it remains unknown whether 
$\IA_n$ is finitely presentable for $n>3$. In 1997, Krsti\'c-McCool 
proved that $\IA_3$ is not finitely presentable \cite[Theorem 1]{km}.

Magnus's set of generators for $\IA_n$ is the collection of automorphisms 
\[
\{g_{ij}\co 1 \leq i,j \leq n, i \neq j\} \cup 
\{ f_{ijk}\co 1 \leq i,j,k \leq n, j \neq i \neq k, i<j \} \]
where $g_{ij}$ is defined by 
\[ 
g_{ij}(x_r) = 
\left \{
\begin{array}{ll}
x_r & r \neq j \\
x_i x_j x_i^{-1} & r = j
\end{array}
\right.
\]
and $f_{ijk}$ is defined by 
\[ 
f_{ijk}(x_r) = 
\left \{
\begin{array}{ll}
x_r & r \neq k \\
x_k [x_i, x_j] & r = k.
\end{array}
\right.
\]
Recently Cohen-Pakianathan \cite{cp},
Farb \cite{farb}, and Kawazumi \cite{kawa} independently determined the 
first homology and cohomology groups of $\IA_n$ using a map based on 
Johnson's homomorphism for the Torelli group \cite{johnhom}. 
For a group $G$, let $G^{(1)} = G$, and 
let $G^{(k)}$ be $k$th term of the lower central series. 
Let $H_1(G)$ denote the integral homology of $G$, and $H_1(G)^*$ its dual.
Cohen-Pakianathan \cite[Corollary 3.5]{cp}, 
Farb \cite[Theorem 1.1]{farb}, and Kawazumi \cite[Theorem 6.1]{kawa} 
proved the following: 
\begin{thm}\label{jf}
There exists a surjective homomorphism 
\[ \tau\co {\rm IA}_n \ra \Wedge^2 H_1(F_n) \otimes H_1(F_n)^*. \]
Furthermore, the induced map 
\[ \tau_*\co H_1({\rm IA}_n) \ra \Wedge^2 H_1(F_n) \otimes
H_1(F_n)^* \]
is an isomorphism of {\rm GL}$(n,\mathbb{Z})$-modules. 
In particular, the free abelian group $H_1${\rm (IA}$_n${\rm )} has rank 
\[ {\rm rank}(H_1({\rm IA}_n)) = n^2(n-1)/2\,. \]
\end{thm}
\noindent See Section \ref{higherdef} for the definition of $\tau$. Note 
that the theorem implies that $\IA_n^{(2)}$ is equal to the kernel of $\tau$. 
Andreadakis \cite{andrea} proved Theorem \ref{jf} for the case $n=3$.

Recall that $H_1(F_n) \simeq F_n/F_n^{(2)} \simeq \mathbb{Z}^n$ has 
a standard basis 
$\{e_1,\ldots,e_n\}$ with dual basis $\{e_1^*,\ldots,e_n^*\}$.
The image of the Magnus generators under $\tau$ is 
\[ \tau(g_{ij}) = (e_i \wedge e_j) \otimes e_j^* \]
$$ \tau(f_{ijk}) = (e_i \wedge e_j) \otimes e_k^*. \leqno{\hbox{and}}$$
The subgroup $\IA_n$ of ${\rm Aut}(F_n)$ is a free group analogue of the 
Torelli subgroup of the based mapping class group. 
The corresponding 
subgroup $\OA_n$ of the outer automorphism group Out$(F_n)$, 
gives a free group analogue of the unbased Torelli group. 
The group of inner automorphisms of $F_n$ is a normal 
subgroup of $\IA_n$, so there is an exact sequence 
\[
1 \ra {\rm OA}_n \ra {\rm Out}(F_n) \ra 
{\rm GL}(n,\mathbb{Z}) \ra 1 
\]
and a surjective homomorphism $\overline{\tau}$ from $\OA_n$ to a quotient of 
$\Wedge^2 H \otimes H^*$ 
\[ \overline{\tau}\co {\rm OA}_n \ra \big(\Wedge^2 H_1(F_n) \otimes H_1(F_n)^* \big) /
H_1(F_n). \]
Almost nothing is known about higher cohomology of $\IA_n$; it is not even 
known whether 
$H^2(\IA_n$) is of finite or infinite 
dimension. 
Our main goal here is to determine the second rational cohomology of 
$\IA_n$ and $\OA_n$ coming from the homomorphisms $\tau$. 
In papers by Hain \cite{hain} and Sakasai \cite{sak}, the Johnson
homomorphism is used to 
understand the second and third rational cohomology groups of the Torelli group. 
In this paper, we use their method to compute the kernel of the
homomophism 
\[
\tau^*\co 
H^2(\Wedge^2 H_1(F_n) \otimes H_1(F_n)^*,\mathbb{Q}) \ra 
H^2({\rm IA}_n,\mathbb{Q})
\]
and an analogue for $\OA_n$.
\begin{thm}\label{tau2} 
The kernel of the homomorphism 
\[
\tau^*\co H^2(\Wedge^2 H_1(F_n) \otimes H_1(F_n)^*, \mathbb{Q}) 
\ra H^2({\rm IA}_n,\mathbb{Q}) 
\]
decomposes into simple 
{\rm GL(}$n,\mathbb{Q}{\rm )}$-modules as 
\[ {\rm ker}(\tau^*) \simeq 
\Phi_{0,\ldots,0,1,0,-1} \oplus \Phi_{1,0,\ldots,0,1,1,-2}. \]
The kernel of the homomorphism
\[
\overline{\tau}^*\co H^2((\Wedge^2 H_1(F_n) \otimes H_1(F_n)^*)/H_1(F_n), \mathbb{Q}) 
\ra H^2({\rm OA}_n,\mathbb{Q})
\] is the simple module 
\[
{\rm ker}(\overline{\tau}^*) \simeq \Phi_{1,0,\ldots,0,1,1,-2}.
\]
\end{thm}
\noindent The notation requires some explanation. 
Briefly, each $\Phi_\alpha$ represents an isomorphism class of a 
simple GL$(n,\mathbb{Q})$-module; the class is 
uniquely determined by the subscript $\alpha$. This notation 
is described in more detail in Section \ref{notation} below 
and is explained comprehensively in Fulton-Harris \cite{fh}. 

A key input is a ``higher'' Johnson homomorphism $\tau^{(2)}$, 
defined on ker($\tau$): 
\begin{thm}\label{higher} There exists a homomorphism 
\[ \tau^{(2)}\co {\rm IA}_n^{(2)} \ra 
\big( \Wedge^2 H_1(F_n) \otimes H_1(F_n)/ \Wedge^3 H_1(F_n) \big) \otimes 
H_1(F_n)^*\]
and a homomorphism 
\[ \overline{\tau}^{(2)}\co {\rm OA}_n^{(2)} \ra 
\big( \big( 
\Wedge^2 H_1(F_n) \otimes H_1(F_n)/ \Wedge^3 H_1(F_n) \big) \otimes 
H_1(F_n)^* \big) / \Wedge^2 H_1(F_n). \]
The kernel of $\tau^{(2)}$ (respectively, $\overline{\tau}^{(2)}$) 
contains {\rm IA}$_n^{(3)}$ (respectively, {\rm OA}$_n^{(3)}$) 
as a subgroup of finite index. 
Tensoring each target space with $\mathbb{Q}$, 
the image of {\rm IA}$_n$ in 
\[
\big( \Wedge^2 H_1(F_n) \otimes H_1(F_n)/ \Wedge^3 H_1(F_n) \big) \otimes 
H_1(F_n)^* \otimes \mathbb{Q}
\] is isomorphic as a {\rm GL(}$n,\mathbb{Q}${\rm )}-representation to 
\[
\Phi_{0,1,0,\ldots,0} \oplus \Phi_{1,1,0,\ldots,0,1,-1}
\]
and the image of {\rm OA}$_n$ in 
\[
\big( \big( \Wedge^2 H_1(F_n) \otimes H_1(F_n)/ \Wedge^3 H_1(F_n) \big) \otimes 
H_1(F_n)^* / \Wedge^2 H_1(F_n) \big) \otimes \mathbb{Q}
\] is isomorphic to $\Phi_{1,1,0,\ldots,0,1,-1}$.\end{thm}

Theorem \ref{higher} is an analogue of a result of Morita's 
for the Torelli group \cite{morita}. Satoh has determined the 
image of $\tau^{(2)}$ over $\mathbb{Z}$ in recent work \cite{sat}. 

An interesting filtration $K_n^{(i)}$ of $\IA_n$ is given by the subgroups 
\[
K_n^{(i)} = {\rm ker(Aut}(F_n) \ra {\rm Aut}(F_n/F_n^{(i)})) 
\]
Following from the definitions and Theorem \ref{jf}, we have 
$K_n^{(1)}=\IA_n$ and $K_n^{(2)} = [\IA_n,\IA_n]$. 
Andreadakis proved that $K_n^{(i)}$ is 
equal to the lower central series in the case $n=2$ 
\cite[Theorem 6.1]{andrea}, and 
that the first three terms of these series coincide in the case $n=3$ 
\cite[Theorem 6.2]{andrea}. This 
led him to conjecture that these two series might be equal for all $n$. 
We reprove his result for $n=3$ and show that the quotient 
$K_n^{(3)}/\IA_n^{(3)}$ is finite for all $n$. 
\begin{cor}\label{rank}
The quotient $K_n^{(3)}/{\rm IA}_n^{(3)}$ is finite for all 
$n$ and trivial for $n=3$. 
The rank of the free abelian group $K_n^{(2)}/K_n^{(3)}$ is 
\[
{\rm rank}(K_n^{(2)}/K_n^{(3)}) = 
\frac{1}{3} n^2(n^2-4) + \frac{1}{2} n(n-1).
\]
\end{cor}
Andreadakis computed the rank of $K_n^{(2)}/K_n^{(3)}$ 
only for the case $n=3$, where he determined the rank to be $18$ 
\cite[Theorem 6.2]{andrea}.

The author would like to thank her advisor Professor Benson Farb for his 
enthusiasm, encouragement, and 
direction. She is also grateful to Fred Cohen, Daniel Biss, David Ogilvie, and 
Ben Wieland for their helpful explanations and suggestions, and to 
Professor Richard Hain for explaining how to detect elements in the 
kernel of $\tau^*$. The author was partially supported 
by an NSERC PGSA fellowship.


\section{Preliminaries}

\subsection{Notation}\label{notation}
The reader is referred to Fulton-Harris \cite{fh} for the representation theory
in this paper. Of particular relevance are Sections 1, 15 and 24.

Any finite dimensional representation of 
GL($n,\mathbb{C}$) decomposes as a tensor product of 
a representation of SL($n,\mathbb{C}$) with a one-dimensional representation. 
Let $\mathfrak{gl}(n,\mathbb{C})$ (respectively, $\mathfrak{sl}(n,\mathbb{C})$) 
be the Lie algebra of GL($n,\mathbb{C}$) (respectively, SL($n,\mathbb{C}$)). 
The finite-dimensional representations of 
SL($n,\mathbb{C}$) and $\mathfrak{sl}(n,\mathbb{C})$ are in natural bijection. 
Let $L_i$ be the element in the dual of the Cartan subalgebra 
of diagonal matrices in $\mathfrak{sl}(n,\mathbb{C})$ which sends a matrix
to its $i$th diagonal entry. Any irreducible representation of 
GL$(n,\mathbb{C})$ is determined up to isomorphism by a 
unique highest weight for SL$(n,\mathbb{C})$ and a power of the 
determinant. This information corresponds to a sum 
\[ a_1L_1 + a_2(L_1 + L_2) + \cdots + a_n(L_1 + L_2 + \cdots + L_n) \] 
where $a_i$ is a non-negative integer, $1 \leq i \leq n-1$, and 
$a_n \in \mathbb{Z}$. We let $\Phi_{a_1,a_2,\ldots,a_n}$ 
denote the representation corresponding to this sum. 
For example, the standard representation of GL$(n,\mathbb{C})$ is 
given by $\Phi_{1,0,\ldots,0}$, and its dual is given by 
$\Phi_{0,\ldots,0,1,-1}$. 

All of these representations are defined over $\mathbb{Q}$, so 
they can be considered as irreducible representations of 
GL($n,\mathbb{Q}$) and $\mathfrak{gl}(n,\mathbb{Q})$.


\subsection{Deriving an exact sequence}\label{seq}
The following exact sequence will be used in the computations underlying 
Theorem \ref{tau2}. 
While this sequence is implicit in Sullivan's paper \cite{sul}, 
we include its derivation here for completeness.
\begin{lem}\label{sequence}
For any group $G$ with finite rank $H^1(G,\mathbb{Q})$, 
the following sequence is exact: 
\[
\xymatrix{
0 \ar[r] & {\rm Hom}(G^{(2)}/G^{(3)},\mathbb{Q}) \ar[r] &
\Wedge^2 H^1(G,\mathbb{Q}) 
\ar[r]^{\cup} & H^2(G,\mathbb{Q}) 
}
\]
\end{lem}
\noindent 
The left hand map is induced by the surjective map 
\begin{align*}
\Wedge^2 G_{ab}& \ra G^{(2)}/G^{(3)}\\
g_1 \wedge g_2&\mapsto [g_1,g_2] G^{(3)}
\end{align*} and the isomorphism 
\[
\Wedge^2 H^1(G,\mathbb{Q}) \simeq {\rm Hom}(\Wedge^2 G_{ab},\mathbb{Q}).
\] The right hand map of the sequence comes from the cup product on cohomology.

Given an exact sequence of groups 
\[
1 \ra G^{(2)} \ra G \ra G_{ab} \ra 0 
\] the Hochschild-Serre spectral sequence yields the five-term exact
sequence 
\begin{align*}
0 \ra 
H^1(G_{ab},H^0(G^{(2)},\mathbb{Q})) \ra 
H^1(G,\mathbb{Q})& \ra H^0(G_{ab},H^1(G^{(2)},\mathbb{Q}) ) \ra\\
\ra\,
&H^2(G_{ab},H^0(G^{(2)},\mathbb{Q})) \ra H^2(G,\mathbb{Q}).
\end{align*}
(For details, consult Section 7.6 of Brown's book \cite{brown}, 
particularly Corollary 6.4.)
Now $H^1(G_{ab},\mathbb{Q})$ and $H^1(G,\mathbb{Q})$ are isomorphic 
and finite rank, so the sequence immediately reduces to 
\[0 \ra H^0(G_{ab},H^1(G^{(2)},\mathbb{Q})) 
\ra H^2(G_{ab},\mathbb{Q}) 
\ra H^2(G,\mathbb{Q}).\] 
We also have $H^2(G_{ab},\mathbb{Q}) = \Wedge^2 H^1(G,\mathbb{Q})$, 
so we are done once we interpret $H^0(G_{ab},H^1(G^{(2)},\mathbb{Q}))$.
By definition, it 
is the set of invariants of $H^1(G^{(2)},\mathbb{Q})) = {\rm
Hom}(G^{(2)},\mathbb{Q})$ under the action of $G_{ab}$. 
This action is given by lifting $G_{ab}$ to $G$ and acting by conjugation.
Thus $H^0(G_{ab},H^1(G^{(2)},\mathbb{Q}))$ 
is the set of homomorphisms of ${\rm Hom}(G^{(2)},\mathbb{Q})$ which
are trivial on $G^{(3)}$, i.e. 
\[H^0(G_{ab},H^1(G^{(2)},\mathbb{Q})) = 
{\rm Hom}(G^{(2)}/G^{(3)},\mathbb{Q}). \]


\subsection{A higher Johnson homomorphism}\label{higherdef}

Unless otherwise stated, $n \geq 3$ will be fixed, IA will denote $\IA_n$, 
and $F$ will denote $F_n$. 

In this subsection, we define a map 
\[ \tau^{(2)}\co {\rm IA}^{(2)} \ra 
\big( \Wedge^2 H_1(F) \otimes H_1(F)/ \Wedge^3 H_1(F) \big) \otimes 
H_1(F)^* \]
which is analogous to the Johnson homomorphism for IA. A generalization of 
this map to all $\IA^{(n)}$ and to some related classes of groups is discussed 
in \cite{cp}, \cite{kawa}, and \cite{sat}.

We begin with a map 
\[
\delta^{(2)}\co {\rm IA}^{(2)} \ra {\rm Hom}(H_1(F),F^{(3)}/F^{(4)}) 
\]
defined as follows: 

For $f \in {\rm IA}^{(2)}$, map $f$ to 
$(x \mapsto f(\tilde{x}) \tilde{x}^{-1})$
where $\tilde{x}$ is a lift of $x \in F/F^{(2)}$ to $F/F^{(4)}$. 
From the fact that $\IA^{(2)}$ acts as the identity on both 
$F/F^{(2)}$ and $F^{(2)}/F^{(4)}$, it is straightfoward 
to verify that $\delta^{(2)}(f)$ is a well-defined homomorphism 
from $H_1(F)$ to $F^{(2)}/F^{(4)}$. 
Furthermore, since $F/F^{(4)}$ factors through $F/F^{(3)}$, and since 
$\IA^{(2)}$ also acts as the identity on $F/F^{(3)}$, it follows 
that $f(\tilde{x})\tilde{x}^{-1}$ actually lies in 
$F^{(3)}/F^{(4)}$. Therefore we have $\delta^{(2)}(f) \in$ Hom$(H_1(F),F^{(3)}/F^{(4)})$. 

The GL$(n,\Z)$-module Hom($H_1(F),F^{(3)}/F^{(4)}$) is isomorphic to 
\[ \big( \Wedge^2 H_1(F) \otimes H_1(F)/ \Wedge^3 H_1(F) \big) \otimes 
H_1(F)^*. \] We define the map 
$\tau^{(2)}$ to be the composition of $\delta^{(2)}$ with this isomorphism.

Since $\IA^{(3)}$ acts as the identity on $F/F^{(4)}$, it lies in 
the kernel of $\tau^{(2)}$.

\section{The kernel of the Johnson homomorphism on \newline second cohomology}
There are three main steps to the proof of Theorem \ref{tau2}, 
and we will first 
restrict our attention to $\IA$, postponing $\OA$ until subsection \ref{conclusion}. 
First we decompose 
$(H_1(F) \Wedge H_1(F)) \otimes H_1(F)^*$ into simple 
GL($n,\mathbb{Q}$)-submodules. The kernel of $\tau$ is 
a GL($n,\mathbb{Q}$)-submodule, so is a direct sum of a subset of 
these simple submodules. 
Next we use the exact sequence of 
Section \ref{seq} to identify 
isomorphism classes of simple modules contained in the kernel of $\tau^*$.
Finally, we use the method of 
Hain \cite{hain} and Sakasai \cite{sak} to find classes which are not in the
kernel. Steps 2 and 3 
exhaust the set of classes in the decomposition of $(H_1(F) \Wedge
H_1(F)) \otimes H_1(F)^*$, showing that 
the terms identified in the step 2 are the only ones in the kernel. 

Once the computations underlying Theorem \ref{tau2} are complete, 
Theorem \ref{higher} and Corollary \ref{rank} are straightforward. 
The fourth subsection below contains their proofs. 

\subsection{Module decompositions}
We begin with a decomposition of the relevant GL$(n,\mathbb{Q})$-modules
into direct sums of simple submodules. 
Let $H = H_1(F)$ and 
\[
U = \Wedge^2 H \otimes H^*.
\]
Let $H_\mathbb{Q}$ and  $U_\mathbb{Q}$ 
denote the tensor product 
with $\mathbb{Q}$ of the spaces $H$ and $U$, respectively. 
The action of GL$(n,\mathbb{Q})$ on $H_\mathbb{Q}^*$ 
will be dual to that on $H_\mathbb{Q}$.

Direct computation yields the following two lemmas. 
\begin{lem}\label{decomp} $U^*_\mathbb{Q}$ decomposes into a direct sum of
simple {\rm GL}$(n,\mathbb{Q})$-submodules as: 
\[
U^*_\mathbb{Q} \simeq 
\left \{
\begin{array}{ll}
\Phi_{2,0,-1} \oplus \Phi_{0,1,-1} & n=3 \\
\Phi_{1,0,\ldots,0,1,0,-1} \oplus \Phi_{0,\ldots,0,1,-1} & n \geq 3
\end{array}
\right.
\] \end{lem}
\def\Strut{\vrule width 0pt depth 7pt}
\begin{lem}\label{wedge2} As a {\rm GL}$(n,\mathbb{Q})$-module, 
$\Wedge^2U^*_{\mathbb{Q}}$ decomposes as: 
\[ 
\Wedge^2 U^*_{\mathbb{Q}} \simeq 
\left \{
\begin{array}{ll}
\Strut2\Phi_{1,0,-1} \oplus 2\Phi_{2,1,-2} & n=3 \\
\Strut\Phi_{2,0,0,-1} \oplus 3\Phi_{0,1,0,-1} \oplus \Phi_{0,3,0,-2} \oplus
\Phi_{3,0,1,-2} \oplus 2\Phi_{1,1,1,-2} & n=4 \\
2\Phi_{1,1,0,0,-1} \oplus 3\Phi_{0,0,1,0,-1} \oplus \Phi_{0,1,2,0,-2}\, \oplus \\ 
\quad \quad \quad \quad \quad \quad \quad \quad \quad \quad \quad \quad \oplus\,
\Phi_{2,1,0,1,-2}
\oplus 2\Phi_{1,0,1,1,-2} & n=5 
\end{array}
\right.
\] For $n=6$:
\begin{align*}
\Wedge^2 U^*_{\mathbb{Q}} \simeq 
\Phi_{0,2,0,0,0,-1} \oplus 2\Phi_{1,0,1,0,0,-1}& 
\oplus 3\Phi_{0,0,0,1,0,-1} \oplus\\
&\oplus \Phi_{0,1,0,2,0,-2} \oplus \Phi_{2,0,1,0,1,-2} \oplus 2\Phi_{1,0,0,1,1,-2} 
\end{align*}
For $n \geq 7$:
\begin{align*}
\Wedge^2 U_\mathbb{Q}^* \simeq 
\Phi_{0,1,0,\ldots,0,1,0,0,0,-1}& \oplus 
2\Phi_{1,0,\ldots,0,1,0,0,-1} \oplus
3\Phi_{0,\ldots,0,1,0,-1} \oplus\\
&\oplus \Phi_{0,1,0,\ldots,0,2,0,-2} \oplus \Phi_{2,0,\ldots,0,1,0,1,-2} \oplus 
2\Phi_{1,0,\ldots,0,1,1,-2}
\end{align*}
\end{lem}


\subsection{A lower bound on ker${\bf (\tau^*)}$}\label{lower}
The key to finding terms in the kernel of $\tau^*$ comes from the
following diagram:
\[
\xymatrix{  & & \Wedge^2 H^1(U,\mathbb{Q}) \ar[d]^{\tau^*} \ar[r]^{\cup} & 
H^2(U,\mathbb{Q})  \ar[d]^{\tau^*} \\
0 \ar[r] & {\rm Hom}({\rm IA}^{(2)}/{\rm IA}^{(3)},\mathbb{Q}) \ar[r] & 
\Wedge^2 H^1({\rm IA},\mathbb{Q}) \ar[r]^{\cup} & H^2({\rm IA},\mathbb{Q}) }
\] The bottom line is the sequence of Lemma \ref{sequence}. The right hand
horizontal maps are the 
cup product on cohomology; the top one is an isomorphism since $U$ is
abelian. The diagram commutes by naturality of the cup product. The
vertical maps are induced by 
$\tau$ on second cohomology. 
The left vertical map is an isomorphism by Theorem \ref{jf}, 
so there exists an
injective homomorphism from Hom$({\rm IA}^{(2)}/{\rm IA}^{(3)},\mathbb{Q})$ to
$\Wedge^2H^1(U,\mathbb{Q})$. Every map is equivariant with respect to the action of 
${\rm Aut}(F)$. We can deduce from the diagram that the kernel of 
\[
\tau^*\co H^2(U,\mathbb{Q}) \ra H^2({\rm IA},\mathbb{Q})
\] is precisely 
\[
{\rm ker}(\tau^*) = {\rm image \big(Hom}({\rm IA}^{(2)}/{\rm IA}^{(3)},\mathbb{Q}) \big)
\subset \Wedge^2 H^1(U,\mathbb{Q})
\] a GL$(n,\mathbb{Q})$-submodule of $\Wedge^2 H^1(U,\mathbb{Q})$.

A first approximation of this image is achieved by considering the map $\tau^{(2)}$ 
defined in Section \ref{higherdef}. The target space 
Hom$(H,F^{(2)}/F^{(3)}) \otimes \mathbb{Q}$ is isomorphic as a 
GL$(n,\mathbb{Q})$-module to 
\[ 
\big( \Wedge^2 H_\mathbb{Q} \otimes H_\mathbb{Q} / \Wedge^3 H_\mathbb{Q} \big) \otimes H_\mathbb{Q}^*
\] 
which decomposes as 
\[ \Phi_{0,1,0,\ldots,0} \oplus \Phi_{1,1,0,\ldots,0,1,-1} \oplus
\Phi_{2,0,\ldots,0}. \]
Consider the element $[g_{12},g_{21}]$ in $\IA^{(2)}$. This acts on the generators of 
$F$ by 
\[ 
x_r \mapsto 
\left \{
\begin{array}{ll}
\Strut[[x_1,x_2],x_r]x_r & r \in \{1,2\} \\
x_r & r \notin \{1,2\}.
\end{array}
\right.
\]
Thus $[g_{12},g_{21}]$ maps by $\tau^{(2)}$ to (the equivalence class of)
\[(e_1\wedge e_2) \otimes e_1 \otimes e_1^* + (e_1\wedge e_2) \otimes e_2
\otimes e_2^* \]
in $\big(\Wedge^2 H_\mathbb{Q} \otimes H_\Q/\Wedge^3 H_\mathbb{Q}\big) \otimes H_\mathbb{Q}^*$. 

Now consider the GL($n,\mathbb{Q}$)-equivariant map 
\[
\big(\Wedge^2 H_\mathbb{Q} \otimes H_\mathbb{Q} / \Wedge^3 H_\mathbb{Q} \big)
\otimes H_\mathbb{Q}^* \ra H_\mathbb{Q}^{\otimes 3} \otimes H_\mathbb{Q}^* 
\] defined by 
\begin{align*}
(a \wedge b) \otimes c \otimes d^* &\mapsto\\
a \otimes b& \otimes c \otimes d^* - b \otimes a \otimes c \otimes d^* + 
c \otimes b \otimes a \otimes d^* - b \otimes c \otimes a \otimes d^*
\end{align*} 
and the GL($n,\mathbb{Q}$)-equivariant map 
\[
H_\mathbb{Q}^{\otimes 3} \otimes H_\mathbb{Q}^* \ra 
\Wedge^2 H_\mathbb{Q}
\] 
$$
a \otimes b \otimes c \otimes d^* \mapsto 
d^*(c) a \wedge b.\leqno{\hbox{defined by}} 
$$
Applying the composition of these maps to 
\[(e_1\wedge e_2) \otimes e_1 \otimes e_1^* + (e_1\wedge e_2) \otimes e_2
\otimes e_2^* \] 
$$
6 (e_1 \wedge e_2)\leqno{\hbox{gives}} 
$$
which is a highest weight vector for the module $\Phi_{0,1,0,\ldots,0,0}$. 

Let $E_{ij}$ be the matrix in $\mathfrak{gl}(n,\mathbb{Q})$ with a $1$ in
the $ij$th place 
and $0$'s everywhere else. 
Applying to the vector 
\[(e_1\wedge e_2) \otimes e_1 \otimes e_1^* + (e_1\wedge e_2) \otimes e_2
\otimes e_2^* \] 
the map
\begin{align*}
(a \wedge b) \otimes c \otimes d^* &\mapsto\\
a \otimes b& \otimes c \otimes d^* - b \otimes a \otimes c \otimes d^* + 
c \otimes b \otimes a \otimes d^* - b \otimes c \otimes a \otimes d^*
\end{align*}
followed by $E_{1n}$, and then by 
\[ a \otimes b \otimes c \otimes d^* \mapsto 
(a \wedge b) \otimes c \otimes d^*
\]
$$-4(e_1 \wedge e_2) \otimes e_1 \otimes e_n^* \leqno{\hbox{we obtain}}$$
This is a highest weight vector for the irreducible representation 
corresponding to $\Phi_{1,1,0,\ldots,0,1,-1}$. 

A simple GL$(n,\mathbb{Q})$-module is equal to the orbit of its highest weight 
by the action of GL$(n,\mathbb{Q})$. 
Since we have found highest weights for $\Phi_{0,1,0,\ldots,0,0}$ and 
$\Phi_{1,1,0,...,0,1,-1}$ 
in the image of the GL$(n,\mathbb{Q})$-invariant map $\tau^{(2)}$, 
both modules are contained in the image of $\tau^{(2)}$. 
Therefore the dual module 
\[
\Phi_{0,\ldots,0,1,0,-1} \oplus \Phi_{1,0,\ldots,0,1,1,-2} \]
is contained in the kernel of $\tau^{*}$. 

\subsection{An upper bound on ker$({\bf \tau^*})$}\label{upper}
In this subsection, we make use of the fact that 
elements in the image of $\tau_*\co H_2({\rm IA},\mathbb{Q}) \ra
H_2(U,\mathbb{Q})$ 
correspond to elements in the dual $H^2({\rm IA},\mathbb{Q})$ which lie outside of
the kernel of $\tau^*$. This will give us an upper bound on ker$(\tau^*)$.

If two elements $g_1$ and $g_2$ of a group $G$ commute, there is a 
homomorphism 
\[
\zeta\co \mathbb{Z}^2 \ra G
\]
defined by mapping the two standard generators of $\mathbb{Z}^2$ to 
$g_1$ and $g_2$, inducing a map 
\[
\mathbb{Z} \simeq H_2(\mathbb{Z}^2) \ra H_2(G).
\] 
The image of the generator $1 \in H_2(\mathbb{Z}^2)$ in $H_2(G)$ is known as an 
{\it abelian cycle}, for which we write 
\[
\{g_1,g_2\} := {\rm image}(\zeta_*(1)) \in H_2(G).
\] 
In general, this element may be zero. In our
case, we identify some abelian cycles in $H_2({\rm IA},\mathbb{Q})$ and compute
their image in $H_2(U,\mathbb{Q})$. These turn out to be nonzero, giving us
nonzero cycles in the image of $\tau_*$ in $H_2({\rm IA},\mathbb{Q})$. 
We list below a set of GL$(n,\mathbb{Q})$-equivariant homomorphisms
which, pre-composed with $\tau$, will detect these cycles. 

First define GL$(n,\mathbb{Q})$-homomorphisms:
\begin{itemize}
\item $f_1\co \Wedge^2(\Wedge^2 H_\mathbb{Q} \otimes H_\mathbb{Q}^*) \ra 
(\Wedge^2 H_\mathbb{Q} \otimes H_\mathbb{Q}^*)^{\otimes 2} $ \newline
$
(a_1 \wedge b_1 \otimes c_1^*) \wedge (a_2 \wedge b_2 \otimes c_2^*)
\mapsto$ 
\newline
\mbox{
$\quad \quad
((a_1 \wedge b_1 \otimes c_1^*) \otimes (a_2 \wedge b_2 \otimes c_2^*)) - 
((a_2 \wedge b_2 \otimes c_2^*) \otimes (a_1 \wedge b_1 \otimes c_1^*))$
}\vspace{0.05in}

\item
$f_2\co (\Wedge^2 H_\mathbb{Q} \otimes H_\mathbb{Q}^*)^{\otimes 2} \ra 
(H_\mathbb{Q}^{\otimes 2} \otimes H_\mathbb{Q}^*)^{\otimes 2}$ \newline
$
(a_1 \wedge b_1 \otimes c_1^*) \otimes (a_2 \wedge b_2 \otimes c_2^*)
\mapsto$ \newline 
\mbox{
$
\quad \quad \quad \quad \quad
((a_1 \otimes b_1 - b_1 \otimes a_1) \otimes c_1^*) \otimes 
((a_2 \otimes b_2 - b_2 \otimes a_2) \otimes c_2^*)) $
}
\end{itemize}
and let $f$ be the composition:
\[
f = f_2 \circ f_1\co \Wedge^2 U_{\mathbb{Q}} \ra (H_\mathbb{Q}^{\otimes 2} \otimes
H_\mathbb{Q}^*)^{\otimes 2}
\] Define a vector $v \in 
(H_\mathbb{Q}^{\otimes 2} \otimes H_\mathbb{Q}^*)^{\otimes 2}$ by 
\[
v= (a_1 \otimes b_1 \otimes c_1^*) \otimes (a_2 \otimes b_2 \otimes c_2^*) . 
\] The following set of 
GL($n,\mathbb{Q}$)-homomorphisms are defined on the image of $f$:
\begin{itemize}
\item 
$g_1(v) = c_1^*(b_2) c_2^*(b_1) (a_1 \wedge a_2) \ \in \ \Wedge^2 H_\mathbb{Q}$
\vspace{0.05in}
\item 
$g_2(v) = c_1^*(b_1) c_2^*(b_2) (a_1 \wedge a_2) \ \in \ \Wedge^2 H_\mathbb{Q}$
\vspace{0.05in}
\item 
$h_1(v) = c_1^*(b_1) (a_1 \otimes a_2 \otimes b_2 \otimes c_2^*) \ \in \ 
H_\mathbb{Q}^{\otimes 3} \otimes H_\mathbb{Q}^*$ 
\vspace{0.05in}
\item 
$h_2(v) = c_1^*(b_2) (a_1 \otimes a_2 \otimes b_1 \otimes c_2^*) \ \in \ 
H_\mathbb{Q}^{\otimes 3} \otimes H_\mathbb{Q}^*$ 
\vspace{0.05in}
\item 
$k(v) = (a_1 \wedge b_1) \otimes (a_2 \wedge b_2) \otimes (c_1^* \wedge
c_2^*) \ \in \ 
(\Wedge^2 H_\mathbb{Q})^{\otimes 2} \otimes \Wedge^2 H_\mathbb{Q}^*$
\vspace{0.05in}
\item 
$l(v) = (a_1 \wedge a_2 \wedge b_1) \otimes b_2 \otimes c_1^* \otimes
c_2^* \ \in \
\Wedge^3 H_\mathbb{Q} \otimes H_\mathbb{Q} \otimes (H_\mathbb{Q}^*)^{\otimes 2}$
\vspace{0.05in}
\item 
$m(v) = c_1^*(b_1)(a_1 \otimes (a_2 \wedge b_2) \otimes c_2^*) \ \in \ 
H_\mathbb{Q} \otimes \Wedge^2 H_\mathbb{Q} \otimes H_\mathbb{Q}^*$
\vspace{0.05in}
\item 
$n(v) = (a_1 \wedge a_2 \wedge b_1 \wedge b_2) \otimes (c_1^* \wedge
c_2^*) \ \in \ 
\Wedge^4 H_\mathbb{Q} \otimes \Wedge^2 H_\mathbb{Q}^*$
\end{itemize}

This completes the list of required homomorphisms. 
We work through proving that there are two terms in the isomorphism 
class of $\Phi_{0,\ldots,0,1,0,-1}$ in some detail. For the remaining terms,
only the abelian cycles and the homomorphisms needed to perform the
computations are given. 

Consider the abelian cycles 
\begin{align*}
\omega_1& = \{g_{12}g_{1n},g_{2n}\}\\
\omega_2& = \{g_{13},g_{2n}\}\\
\omega_3& = \{f_{12n},g_{3(n-1)}\} 
\end{align*} 
of $H_2({\rm IA},\mathbb{Q})$. 
For $\omega_2$ and $\omega_3$, we must have $n>3$.
By direct computation, we have 
\[\tau_*(\omega_1) = 
(e_1 \wedge e_2 \otimes e_2^* + e_1 \wedge e_n \otimes e_n^*) 
\wedge (e_2 \wedge e_n \otimes e_n^*). \] 
Applying $g_1 \circ f$, we obtain $2 (e_1 \wedge e_2)$, a highest weight
for the 
representation $\Phi_{0,1,0,\ldots,0,0}$, implying that the 
dual representation $\Phi_{0,\ldots,0,1,0,-1}$ is not in the kernel of $\tau^*$. 

Now compute:
\[\tau_*(\omega_2) = 
(e_1 \wedge e_3 \otimes e_3^*) \wedge (e_2 \wedge e_n \otimes e_n^*) \] 
Notice that this vector is in the kernel of $g_1 \circ f$. Applying 
$g_2 \circ f$ to $\tau_*(\omega_2)$, we again obtain $2 (e_1 \wedge e_2)$, giving us a 
second submodule of 
$H^2(U,\mathbb{Q})$ in the class $\Phi_{0,\ldots,0,1,0,-1}$
not in ker$(\tau^*)$.

The simple module $\Phi_{1,0,\ldots,0,1,0,0,-1}$ appears with
multiplicity two in $H^2(U,\mathbb{Q})$, so we can use 
\[
h_2 \circ f \circ E_{3n} \circ \tau_*(\omega_1)
\]
$$h_1 \circ f \circ \tau_*(\omega_3)\leqno{\hbox{and}}$$
to show that both of these terms are in the image of $\tau_*$.

We list here four more GL$(n,\mathbb{Q})$-homomorphisms and abelian cycles which give 
highest weights of submodules whose dual modules are not in ker$(\tau^*)$:
\begin{itemize}
\item
$k \circ f \circ E_{2(n-1)} \circ E_{1n} \circ \tau_*(\omega_1)$ 
is a highest weight for the submodule of $H_2(U,\mathbb{Q})$ whose dual is 
$\Phi_{0,1,0,\ldots,0,2,0,-2}$ in $H^2(U,\mathbb{Q})$.
\vspace{0.05in}
\item
$l \circ f \circ E_{3n} \circ E_{1n} \circ \tau_*(\omega_2)$ 
is a highest weight for the module whose dual is $\Phi_{2,0,\ldots,0,1,0,1,-2}$.
\vspace{0.05in}
\item
$m \circ f \circ E_{14} \circ E_{4n} \circ \tau_*(\omega_1)$ 
is a highest weight for the module whose dual is $\Phi_{1,0,\ldots,0,1,1,-2}$.
\vspace{0.05in}
\item
$n \circ f \circ E_{3(n-1)} \circ E_{4n} \circ \tau_*(\omega_2)$ 
is a highest weight for the module whose dual is $\Phi_{0,1,0,\ldots,0,1,0,0,0,-1}$.
\end{itemize}


\subsection{Conclusion of the proofs}\label{conclusion}
The computations needed to prove the theorems for IA are complete. 
Those of the previous subsection show that the submodule 
\begin{align*}
2\Phi_{0,\ldots,0,1,0,-1}& \oplus 2\Phi_{1,0,\ldots,0,1,0,0,-1} \oplus 
\Phi_{0,1,0,\ldots,0,2,0,-2} \oplus \\
&\oplus \Phi_{2,0,\ldots,0,1,0,1,-2} \oplus
\Phi_{1,0,\ldots,0,1,1,-2} \oplus \Phi_{0,1,0,\ldots,0,1,0,0,0,-1}
\end{align*}
of $H^2(U,\mathbb{Q})$ 
is not in the kernel of $\tau^*$. Recall that in subsection \ref{lower},
we proved that the submodule 
\[
\Phi_{0,\ldots,0,1,0,-1} \oplus \Phi_{1,0,\ldots,0,1,1,-2}
\] 
is in the kernel. Comparing this with Lemma \ref{wedge2}, 
this accounts for all terms of the decomposition of $H^2(U,\mathbb{Q})$.
This completes the proof of Theorem \ref{tau2} for IA.

Recall that the map 
\[
\tau^{(2)}\co {\rm IA}^{(2)} \ra {\rm Hom}(H,F^{(3)}/F^{(4)})
\] 
factors through $\IA^{(2)}/{\rm IA}^{(3)}$. The proof of Therem \ref{tau2} 
together with the computations of subsection \ref{lower} 
imply that as $\mathbb{Q}$-vector spaces: 
\[
{\rm IA}^{(2)}/{\rm IA}^{(3)} \otimes \mathbb{Q} \simeq 
\Phi_{0,1,0,\ldots,0} \oplus \Phi_{1,1,0,\ldots,0,1,-1} 
\subset \tau^{(2)}({\rm IA}^{(2)}/{\rm IA}^{(3)}) \otimes \mathbb{Q}
\] Therefore the map 
\[
\tau^{(2)} \otimes id\co {\rm IA}^{(2)}/{\rm IA}^{(3)} \otimes \mathbb{Q} \ra 
{\rm Hom}(H,F^{(3)}/F^{(4)}) \otimes \mathbb{Q}
\] 
is injective. Tensoring with $\mathbb{Q}$ kills only torsion, so this 
completes the proof of Theorem \ref{higher} for IA.

From the definitions, the kernel of $\tau^{(2)}$ is exactly $K_n^{(3)}$, so 
it follows from Theorem \ref{higher} that $K_n^{(3)}$ contains $\IA^{(3)}$ as a finite 
index subgroup. That they are equal when $n=3$ follows from computations to 
come in Section \ref{explicit}. 
As $K_n^{(2)}/K_n^{(3)}$ is a free abelian group, its rank is equal to the dimension of 
the image of $\tau^{(2)} \otimes id$ 
in Hom$(H,F/F^{(3)}) \otimes \mathbb{Q}$, which can be 
computed using the Weyl character formula (see \cite[page 400]{fh}). This concludes 
the proof of Corollary \ref{rank}. 

We now turn our attention to $\OA$, beginning with an analogue of Theorem \ref{jf}. 
Recall that $\OA$ is equal to $\IA/{\rm Inn}(F)$, the quotient of $\IA$ by the 
inner automorphisms of $F$. 
The equivariance of $\tau$ with respect to the action of ${\rm Aut}(F)$ implies that 
${\rm Inn}(F)$ maps to a GL$(n,\mathbb{Z})$-submodule of $\Wedge^2 H \otimes H^*$. 
Thus we obtain an ${\rm Aut}(F)$-equivariant map 
$\overline{\tau}$ defined on $\OA$ having as its image the quotient of 
$\Wedge^2 H \otimes H^*$ by the image of ${\rm Inn}(F)$. 
Note that ${\rm Inn}(F)$ is generated by the elements 
\[
g_{k1} g_{k2} \cdots g_{kn} \quad \quad 1 \leq k \leq n
\] 
and so its image in $\Wedge^2 H \otimes H^*$ is generated by the elements 
\[
\sum_{i=1}^n e_k \wedge e_i \otimes e_i^* \quad \quad 1 \leq k \leq n.
\] 
There is an injective homomorphism 
\[
H \ra \Wedge^2 H \otimes H^*
\] 
given by mapping the generators $e_k$ to 
\[
e_k \mapsto \sum_{i=1}^n e_k \wedge e_i \otimes e_i^* \quad \quad 1 \leq k \leq n.
\] 
So $\overline{\tau}$ is the homomorphism 
\[
\overline{\tau}\co {\rm OA} \ra 
\big( \Wedge^2 H \otimes H^* \big) / H
\] 
with kernel $[\OA,\OA]$.

Similarly there is a ``higher'' Johnson homomorphism for $\OA$. 
The homomorphism $\tau^{(2)}$ maps ${\rm Inn}(F)$ to the submodule of 
$(\Wedge^2 H \otimes H / \Wedge^3 H \big) \otimes H^*$ 
generated by (equivalence classes) of the form 
\[
\sum_{i=1}^n (e_j \wedge e_k \otimes e_i) \otimes e_i^* \] 
There is an injective homomorphism 
\begin{align*}
\Wedge^2 H& \ra 
\big( \Wedge^2 H \otimes H / \Wedge^3 H \big) \otimes H^*\\
e_j \wedge e_k& \mapsto \sum_{i=1}^n (e_j \wedge e_k \otimes e_i) \otimes e_i^* 
\end{align*}
and so there is a homomorphism 
\[
\overline{\tau}^{(2)}\co {\rm OA}^{(2)} \ra 
\big( \big( \Wedge^2 H \otimes H / \Wedge^3 H \big) \otimes H^*
\big) / \Wedge^2 H. 
\] 
Now arguments previously applied to IA complete the proofs of 
Theorems \ref{tau2} and \ref{higher}. 


\subsection{Towards computing the kernel of ${\bf \tau}^*$ on $H^3(U,\mathbb{Q})$.}
\label{explicit}
At this stage, it should be possible to find cocyles in $H^3(\OA)$ using Theorem 
\ref{tau2} and the method of Sakasai \cite{sak}. However 
we are not yet equipped to take this step for IA 
because, although we know the isomorphism classes of terms in the kernel of 
$\tau^*\co H^2(U,\Q) \rightarrow H^2(\IA,\Q)$, 
these classes appear with multiplicity greater than one in the decomposition of 
$\Wedge^2 U^*_\mathbb{Q}$. To follow Sakasai's 
method for computing the kernel in $H^3(U,\mathbb{Q})$ of $\tau^*$, 
we would first need to find explicit highest weights for the simple modules 
in the kernel. We do this for the case $n=3$.

The set of nine Magnus generators of $\IA_3$ gives a set of $36$
generators for the finite rank abelian group $\IA_3^{(2)}/\IA_3^{(3)}$. 
For distinct $i,j,k \in \{1,2,3\}$, the following elements of
$\IA_3^{(2)}/\IA_3^{(3)}$ are equal to $0$: 
\begin{align*}
& [g_{ij},g_{ik}]\\
& [g_{ij},g_{ki}] + [g_{ij}, g_{kj}] \\
& [f_{ijk},g_{ij}] + [f_{ijk},g_{ik}] \\
& [f_{ijk},g_{ji}] + [f_{ijk},g_{jk}] \\
& [f_{ijk},g_{ki}] + [f_{ijk},g_{kj}] + [g_{ij},g_{ji}] +
[g_{ik},g_{ji}] \end{align*} 
where the ``$+$'' notation denotes multiplication in 
$\IA_3^{(2)}/\IA_3^{(3)}$. 
These relations imply that a minimal generating set for
$\IA_3^{(2)}/\IA_3^{(3)}$ has no more than $18$ elements. 
Since $K_3^{(3)}$ is equal to the kernel of $\tau^{(2)}$, 
the formula of Corollary \ref{rank} implies that 
18 is the dimension of the image 
of $\IA_3^{(2)}$ in Hom$(H,F_3/F_3^{(3)}) \otimes \mathbb{Q}$. 
As $\IA_3^{(3)}$ is contained in $K_3^{(3)}$, 
the list of relations must be complete. 
Notice that this implies the 
last claim of Corollary \ref{rank}, that $K_3^{(3)}=\IA_3^{(3)}$, which 
was originally shown by Andreadakis \cite{andrea}.

Now we can identify explicit vectors in the image of
${\rm Hom}(\IA_3^{(2)}/\IA_3^{(3)},\mathbb{Q}$) in 
$\Wedge^2 H^1(U,\mathbb{Q})$. For example, the element 
\[(e_2^* \wedge e_3^* \otimes e_3) \wedge (e_3^* \wedge e_2^* \otimes e_2)
+
(e_2^* \wedge e_3^* \otimes e_1) \wedge (e_1^* \wedge e_3^* \otimes e_3) 
\] 
of $\Wedge^2 H^1(U,\mathbb{Q})$ which sends each of 
\[ (e_2 \wedge e_3 \otimes e_3^*) \wedge (e_3 \wedge e_2 \otimes e_2^*) \] 
$$(e_2 \wedge e_3 \otimes e_1^*) \wedge (e_1 \wedge e_3 \otimes e_3^*) 
\leqno{\hbox{and}}$$  
to $1$ and all other standard basis elements of $\Wedge^2
H^1(U,\mathbb{Q})$ to $0$ is in the image of 
${\rm Hom}(\IA_3^{(2)}/\IA_3^{(3)},\mathbb{Q})$. 


\Addresses\recd

\end{document}

